%
%

\magnification=1200

\font\titfont=cmr10 at 14 pt

\font\headfont=cmr10 at 12 pt



\overfullrule=0in

\def\boxit#1{\hbox{\vrule
 \vtop{%
  \vbox{\hrule\kern 2pt %
     \hbox{\kern 2pt #1\kern 2pt}}%
   \kern 2pt \hrule }%
  \vrule}}

  \def\harr#1#2{\ \smash{\mathop{\hbox to .3in{\rightarrowfill}}\limits^{\scriptstyle#1}_{\scriptstyle#2}}\ }

\def\Cn{{\rm Cone}}

 \def\Ed{{\rm Edge}}
\def\spt{{\rm Spt}}
\def\stab{{\rm Stab}}
\def\as{{\rm ASpan}}

\def\bra#1#2{\langle #1, #2\rangle}
\def\bbf{{\bf F}}
\def\bbj{{\bf J}}

\def\ss{\subset}

\def\oa#1{\overrightarrow #1}
\def\dim{{\rm dim}}
\def\dist{{\rm dist}}

\def\span{{\rm Span\,}}

\def\Sym{{\rm Sym}^2}

\def\pr{{\rm pr}}

\def\arr{\longrightarrow}

\def\rn{\bbr^n}

\def\Int{{\rm Int}}

\def\Symn{{\Sym(\rn)}}

 \def\cd{{\cal C}}

\def\Theorem#1{\medskip\noindent {\bf THEOREM \bf #1.}}
\def\Prop#1{\medskip\noindent {\bf Proposition #1.}}
\def\Cor#1{\medskip\noindent {\bf Corollary #1.}}
\def\Lemma#1{\medskip\noindent {\bf Lemma #1.}}
\def\Remark#1{\medskip\noindent {\bf Remark #1.}}
\def\Note#1{\medskip\noindent {\bf Note #1.}}
\def\Def#1{\medskip\noindent {\bf Definition #1.}}

\def\Ex#1{\medskip\noindent {\bf Example \bf    #1.}}

\def\pf{\medskip\noindent {\bf Proof.}\ }
\def\qed{\hfill  $\vrule width5pt height5pt depth0pt$}

\def\mathqed{  $\vrule width5pt height5pt depth0pt$}

   \def\cc{{\cal C}}     
   \def\cp{{\cal P}}

\def\ch{{\cal H}}   
   
\def\cd{{\cal D}}

\def\cp{{\cal P}}
\def\cf{{\cal F}}

\def\vf{\varphi}

\def\wt{\widetilde}

\def\and{\qquad {\rm and} \qquad}
\def\arr{\longrightarrow}

\def\bbr{{\bf R}}

\def\bbz{{\bf Z}}

\def\a{\alpha}

\def\e{\epsilon}

\def\l{\lambda}

\def\s{\sigma}

\def\D{\Delta}

\def\lloc{L^1_{\rm loc}}

\def\SH{{\rm SH}}

\def\Symn{\Sym(\rn)}
 
\def\USC{{\rm USC}}
\def\fa{{\rm\ \  for\ all\ }}

\def\Fa{{\oa F}}

\def\AA{1}
\def\BB{2}
\def\CC{3}
\def\DD{4}
\def\EE{5}
\def\FF{6}
\def\GGG{7}
\def\HH{8}
\def\II{9}

\def\Fv{F^{\rm visc}}
\def\Fd{F^{\rm dist}}

\def\SH#1{{H}^{#1}}
\def\SHV{\SH {\rm visc}}
\def\SHW{\SH {\rm dist}}
\def\SHI{\SH {\rm class}}

\def \esssup#1{{\rm ess} \!\! \sup_{#1}}
\def \essup#1{{\rm ess}  \sup_{#1}}

\centerline{\titfont THE EQUIVALENCE OF VISCOSITY AND DISTRIBUTIONAL}\medskip
\centerline{\titfont  SUBSOLUTIONS FOR CONVEX  SUBEQUATIONS -- }
\medskip

\centerline{\titfont A STRONG BELLMAN PRINCIPLE }
\medskip

\bigskip

\centerline{\titfont F. Reese Harvey and H. Blaine Lawson, Jr.$^*$}
\vglue .9cm
\smallbreak\footnote{}{ $ {} \sp{ *}{\rm Partially}$  supported by
the N.S.F. } 

\vskip .1in

\centerline{\bf ABSTRACT} \medskip
  \font\abstractfont=cmr10 at 10 pt
 
  {{\parindent= .83in
\narrower\abstractfont \noindent
 There are two  useful ways to extend  nonlinear partial differential
inequalities of second order: one uses viscosity theory and the other uses
the theory of distributions. This paper considers  the convex situation where both 
extensions can be applied.  The main result is that under
 a natural   ``second-order completeness''  hypothesis,
 the two sets of extensons are isomorphic, in a sense that is made precise.
 
}}

\vskip.3in

\centerline{\bf TABLE OF CONTENTS} \bigskip

{{\parindent= .1in\narrower\abstractfont \noindent

\qquad \AA.     Introduction.   
 
 \smallskip

\qquad \BB. Differential Constraints -- Two Approaches.

   \smallskip

\qquad \CC.  The Standard   Intersection Theorem for Closed Convex Sets.  

   \smallskip

\qquad \DD.  Another Intersection Theorem for Closed Convex Sets.  

   \smallskip

\qquad \EE.   The Edge of a Convex Set.  

   \smallskip

\qquad \FF.   Subequations Which are Second-Order Complete.  

   \smallskip

\qquad \GGG.   Characterizing When Stable Means Uniformly Elliptic.  

   \smallskip

\qquad \HH.   The Main Result for Variable Coefficients.  

   \smallskip

\qquad \II.   The Equivalence of  Various Notions of Subharmonicity for Linear Equations.  

   \smallskip

}}

\vskip .1in

{{\parindent= .3in\narrower

\noindent
{\bf Appendices: }\medskip

A.   Review of the Relevant Convex Geometry.
\smallskip

}}

\vfill\eject


\noindent{\headfont \ \AA.  Introduction.}
\medskip

There are two quite distinct approaches to the study of nonlinear partial differential
inequalities of second order:  the viscosity approach and the distributional approach.
The purpose of this paper is to prove that in the situation where both can be applied,
there is a natural   ``second-order completeness''  hypothesis  under which the two approaches are,
in a certain precise sense, isomorphic.  

In either approach one can start quite generally by considering the $C^2$-functions $u$
on a manifold $X$ satisfying a second-order constraint $F$. In local coordinates
this  comes down to requiring that 
$$
\left(u(x), D_xu, D^2_x u      \right) \ \in\ F_x \qquad{\rm for\ each\ } x
\eqno{(\AA.1)}
$$
where  the constraint set $F_x$ at $x$ is a subset of the space of 2-jets 
$\bbj^2 = \bbr\times \rn\times\Symn$ (and $\Symn$ is the space of symmetric
$n\times n$-matrices).  In the viscosity case three additional conditions are required
(see Defs. \BB.2 and \HH.1).  When they are satisfied, $F$ is  called a {\bf subequation}.
In the distributional case, {\sl a priori} only linear inequalities make sense, and so
one must require that each $F_x$ be convex.  Thus our starting point for a possible
isomorphism between the two approaches is a {\bf convex subequation} $F$.

The idea now is to extend the {\sl $F$-subharmonicity condition} (\AA.1) to more general ``functions''.
In the viscosity setting we consider the space $\USC(X)$ of upper semi-continuous
$[-\infty, \infty)$-valued functions on $X$. The standard viscosity definition extends the notion
of $F$-subharmonicity to functions $u\in\USC(X)$ by using $C^2$-test functions  
(see Def. \BB.1).  For a distribution $v\in\cd'(X)$ one defines $v$ to be $F$-subharmonic
by requiring that locally $v$ satisfies all the linear  second-order inequalities deducible from $F$.
We shall denote by $\Fv(X)$ and $\Fd(X)$  these two different spaces of $F$-subharmonics.
(See Section \BB\ for a precise formulation.)

In general there is no way to associate a distribution to an upper semi-continuous function
or vice versa.
Elementary examples, such as the subequation $F$ on $\bbr^2$ defined by the inequality
$\partial^2 u/\partial x_1^2 \geq0$, show that $\Fv(X)$ and $\Fd(X)$ can be quite different.
Indeed, for this subequation any upper semi-continuous $u(x_2)$ lies in $\Fv(X)$ 
and any distribution $v(x_2)$ lies in $\Fd(X)$.  The problem is that this subequation $F$ can 
be ``defined using fewer of the independent variables'', a notion made precise in Section \FF\
for any pure second-order constant coefficient subequation.  For this not to happen
is a form of ``completeness'' for the subequation $F$.

The concept we need is formulated as follows, using the standard 2-jet coordinates  
$(r,p,A)\in \bbj^2$.
Given $x\in X$ we say that $F_x$ is {\bf second-order complete}
if for some $r,p$, the associated pure second-order subequation $F_{x,r,p}$ 
cannot be defined using fewer of the independent variables.
 This definition is more robust than it might seem. If it holds for one fibre $F_{x,r,p}$,
 then it holds for all non-empty fibres $F_{x,r',p'}$ at $x$.

Our main result provides a precise isomorphism between $\Fv(X)$ and $\Fd(X)$
under a mild ``regularity''  assumption which is discussed below and in Section 8.

\Theorem {\AA.1} {\sl
Suppose $F$ is a regular convex subequation on a manifold $X$, and that 
$F$ is second-order complete.  \medskip

\item {(A)}  If $u\in \Fv(X)$, then $u\in\lloc(X)$, and as a distribution $u\in \Fd(X)$.

\medskip

\item {(B)}  If $u\in \Fd(X)$, then $u\in \lloc(X)$, and within the $\lloc$-class  $u$
of point-wise defined functions  there exists a unique upper semi-continuous representative
$U\in \Fv(X)$ given by 
$$
U(x) \ \equiv\  {\rm ess}  \limsup_{y\to x} u(y) \ \equiv \ \lim_{r\searrow 0}  {\rm ess}  \sup_{|y| \leq r} u(y).
\eqno{(\AA.2)}
$$
}

\medskip
 It is shown in Section \HH\  that any  convex subequation which is locally affinely jet-equivalent to a 
 constant coefficient subequation (see [HL$_3$]) is regular.  This covers most of the
 non-linear equations that arise in geometry. 
 
More generally, it is shown in Section \HH \ that if the edge of $F$ (see Def. \EE.1) is a vector 
sub-bundle of $J^2(X)$,   then $F$ is regular.
 In particular, this holds if $\Ed(F_x) = \{0\}$ for all $x\in X$.

Theorem \AA.1  is proved by reducing to the linear case.  However, even having done that, 
issues remain.  In order to clarify both this reduction and the problem with the linear case,
we now restrict our discussion to the case of a constant coefficient subequation 
$F\ss \bbj^2 \equiv \bbr\times\rn\times\Symn$ on an open set $X\ss\rn$.

The dual space $(\bbj^2)^*$ of $\bbj^2$ should be viewed as the space of second-order
partial differential operators with constant coefficients.  In both approaches
$$
u\in F(X) \qquad\Rightarrow \qquad Lu\ \geq\ \l
\eqno{(\AA.3)}
$$
for all $L\in (\bbj^2)^*$ and $\l\in\bbr$ which define a closed half-space
$H(L,\l)$ containing $F$. In fact, the converse
$$
Lu\ \geq\ \l \qquad \forall\ H(L,\l) \supset F
 \qquad\Rightarrow \qquad u\in F(X)
\eqno{(\AA.4)}
$$
is also true for both approaches.  For $\Fv(X)$ this is a triviality from the definition,
while for $\Fd(X)$ this is the Definition \BB.3.

The conditions required for $F$ to be a subequation imply that if 
$F\ss H(L,\l)$ with $Lu\equiv \bra a {D^2u} + \bra b {Du} +cu$, then
$$
a\ \geq \ 0\and c\ \leq\ 0.
\eqno{(\AA.5)}
$$
The two facts (\AA.3) and (\AA.4) reduce the isomorphism problem
$\Fv(X)\cong\Fd(X)$ to an isomorphism problem for the associated linear inequalities
$Lu\geq\l$.

However, even in this linear case there are difficulties.  Examples such as 
$$
\left\{ {\partial^2u \over \partial x_1^2}  \geq0  \right\} 
\ \cap\ 
\left\{ {\partial^2u \over \partial x_2^2}  \geq0  \right\},
\qquad{\rm or}\qquad
\left\{ D^2 u\geq0  \right\}
$$ 
are second-order complete, but no isomorphism is possible for all the 
associated linear subequations.  For example, in both cases
$F\ss H(L,\l)$ with $L\equiv \partial^2/\partial x_1^2$ and $\l=0$.

We deal with this by showing that under the hypothesis of second-order completeness,
the convex subequation $F$ can be expressed as an intersection of half-spaces
$H(L,\l)$ where the  associated linear operators are {\bf uniformly elliptic}, i.e., $a>0$ (positive definite).

This is done in detail as follows.  Each convex subset $F$ has an ``edge'' $\Ed(F)$ which
is the largest vector subspace such that $F+\Ed(F)\ss F$.
If $H=H(L,\l)$ is a closed half-space containing $F$, then any rotation of $H$  in the edge directions, no matter how small, will no longer contain $F$.  If for small rotations of $H$ in the directions
orthogonal to $\Ed(F)$, the condition $F\ss H$ is maintained, we say that $H$ is $F$-{\bf stable}.
This condition on $H=H(L,\l)$ holds for one $\l \ \iff $ it holds for generic $\l$, and so we say 
that the linear operator $L$ is $F$-{\bf stable} or $L\in\stab(F)$ (see Def. \DD.1).

Now we can state our result which provides a successful reduction to the linear case.

\Theorem{\AA.2. (Reduction)} {\sl
Suppose that $F\ss \bbj^2$ is a (proper) convex subequation with constant coefficients.
Then $F$ is the intersection of the $F$-stable half-spaces   which contain $F$,   i.e.,
$$
F \ =\ \bigcap _{  \eqalign{&L\in\stab(F) \cr  &F\ss H(L,\l)\cr}   }  H(L,\l).
\eqno{(\AA.6)}
$$

Moreover, if $F$ is second-order complete, then each $F$-stable linear operator $L$ is uniformly
elliptic.  Consequently,
\medskip

\noindent
{\bf (Viscosity)}: Given $u\in\USC(X)$
$$
u\in\Fv(X) 
\qquad\iff\qquad 
Lu \geq_{\rm visc}  \l \quad \forall\ L\in\stab(F) \ \ {\rm with} \ \ H(L,\l)\supset F.
$$

\medskip

\noindent
{\bf (Distributional)}: Given $u\in\cd'(X)$
$$
u\in\Fd(X) 
\qquad\iff\qquad 
Lu \geq_{\rm dist}  \l \quad \forall\ L\in\stab(F)  \ \ {\rm with} \ \ H(L,\l)\supset F.
$$
}

\pf This theorem combines Theorem \DD.2 and Theorem \GGG.2 using the elementary
Remark \CC.3. \qed

\vskip .3in
\centerline
{\bf The Linear Case.}
\medskip

Theorem \AA.2 does not quite reduce the isomorphism problem
 for $\Fv(X)\cong\Fd(X)$ to an isomorphism
problem for the linear case
 $H^{\rm visc}(X) \cong H^{\rm dist}(X)$ where $H=H(L,\l)$ and $L$ is $F$-stable.
Even if we knew  $H^{\rm visc}(X) \cong H^{\rm dist}(X)$ for all $H=H(L,\l)\supset F$ where  
$L$ is $F$-stable (and hence uniformly elliptic) a problem would remain.  Namely,
given $u \in \Fd(X)$, so that $u\in  H^{\rm dist}(X)$ for each such $H$, the associated
upper semi-continuous functions $v_H\in  H^{\rm visc}(X)$ must all be equal, in order
to produce a function $v\in \Fv(X)$.  

This is done in Section \II  \  using a third more classical definition of $H$-subharmonicity
as a bridge between  $H^{\rm dist}(X)$  and $H^{\rm visc}(X)$.
We say $u\in \SHI(X)$ if $u$ is ``sub'' the $H$-harmonics (see Def. \II.1).

\Theorem {\II.3(B)}
{\sl
If  $u\in \SHW(X)$,  then $u\in \lloc(X)$, and within the $\lloc$-class $u$ of point-wise
defined functions, there exists 
 a unique upper semi-continuous representative    $U \in \SHI(X)$.  It is given by}
$$
U(x) \ \equiv \  \overrightarrow{ {\rm ess} \lim_{y\to x} }  u(y) \ \equiv \ \lim_{r\searrow0} {\rm  ess} \!\! \sup_{B_r(x)} u
$$
\medskip

Combined with an actual equality between $H^{\rm class}(X)$  and $H^{\rm visc}(X)$ (Theorem \II.2),
this gives the desired independence, since   the essential-lim-sup-regularization $U$
of $u\in \lloc(X)$  does not depend on  $L$. These results are established in Section \II \
using, among other things,
classical results of Herv\'e-Herv\'e [HH].       With $\SHI(X)$ as a bridge this completes the
proof of Theorem \AA.1 when $F=H=H(L,\l)$ is linear (Cor. \II.4).

\medskip
\noindent
{\bf Some Historical Remarks on the Linear Case.}
The equivalence of $\SHV(X)$ and $\SHW(X)$ for linear elliptic operators 
has been addressed by Ishii [I],
who proves the result for continuous functions but leaves open the case where 
$u\in\SHV(X)$ is a general upper semi-continuous function and the case where
$u\in \SHW(X)$ is a general distribution.
The proof that ``classical implies distributional'' appears in [HH] where the result is proved for 
even more general linear hypoelliptic operators $L$. Other arguments that ``viscosity implies distributional'' 
are known to Hitoshi Ishii and to Andrzej Swiech.
A good discussion of the Greens kernel  appears in ([G]).

\medskip
\noindent
{\bf Proof of Theorem \AA.1.}  The linear case of Theorem \AA.1 combined with the reduction Theorem
\AA.2 yields Theorem \AA.1.

\vfill\eject


\noindent{\headfont \ \BB.  Differential Constraints -- Two Approaches.}
\medskip

Any subset $F\ss \bbj^2 \equiv \bbr\times \rn\times \Symn$
 imposes an unambiguous set of constraints on
the $C^2$-functions on an open subset $X\ss\rn$.  
Given $u\in C^2(X)$ we say that $u$ is $F$-{\bf constrained}  if
$$
\left (u(x), D_x u, D^2_xu \right ) \  \in \ F \fa x\in X,
\eqno{(\BB.1)}
$$
and we let $C^2_F(X)$ denote this set of $C^2$-functions.
If $F$ is a convex set, then $F$ will be referred to as a {\sl convex
constraint}.  A {\sl pure second-order constraint} is a constraint of the
form $F=\bbr\times \rn\times F'$, and in this case (\BB.1) can be 
written more simply as $D^2_xu \in F', \forall \,x\in X$.
We always assume that $F$ is closed.  This ensures that $\cc_F(X)$
is closed in $C^2(X)$.

In and of itself, the condition (\BB.1)  is not  particularly interesting
without expanding the notion to more general ``functions''. 
There are two standard ways of doing this.  Not surprisingly, both 
require (different) additional conditions on the constraining set $F$.
We label these two approaches the ``viscosity'' approach and the 
``distributional'' approach.  Both have their advantages and limitations.
One disadvantage of the distributional approach is that it only makes sense
when $F$ is convex.  A big advantage of the viscosity approach is that convexity
is not required.  However, the ``positivity'' condition described below must be assumed.
In this paper we examine constraint sets $F$ where both approaches apply and characterize
when they are equivalent -- in a sense to be made precise.
\vskip .3in

\centerline{\bf The Viscosity Approach.}
\medskip

The more general ``functions'' considered here are indeed pointwise-defined
functions.  Namely, let $\USC(X)$ denote the space of upper semi-continuous 
$[-\infty, \infty)$-valued functions on $X$.

\medskip
\noindent
{\bf Test Functions:} Given $u\in\USC(X)$ and a point $x\in X$,
a $C^2$-function $\vf$ is a {\sl test function for $u$ at $x$} if $u\leq \vf$
 near $x$ and $u(x)=\vf(x)$. 
 
 \Def{\BB.1} A function $u\in\USC(X)$ is said to be  {\bf $F$-subharmonic on $X$}
 if for each $x\in X$ and each test function $\vf$ for $u$ at $x$, we have 
 $$
\left (\vf(x), D_x \vf, D^2_x \vf  \right ) \ \in \  F 
 $$
The set of all such functions is denoted by $\Fv(X)$.
\medskip

Conditions on the set $F$ are required in order for this definition
to be of any value.  Note that if $\vf$ is a test function for $u$ at $x$,
then so is $\vf(y) + \bra {P(y-x)}{y-x}$ for any $P\in\Symn$ with  $P\geq0$.
This viscosity notion of ``generalized second derivative'' yields a set of possibilities,
which is closed under addition of any $P\geq0$.  Therefore, one must require
the following condition:

\medskip
\noindent
{\bf (Positivity)} 
\centerline
{
$(r,p,A) \ \in\ F \qquad\Rightarrow\qquad (r,p, A+P)\ \in\ F \fa P\geq0.$
\qquad\qquad\qquad
}
\medskip
\noindent
This condition is paramount. In particular, it is both necessary 
and sufficient to ensure that the $C^2$-functions that
are $F$-constrained are in $\Fv(X)$.

We also require the:
\medskip
\noindent
{\bf (Topological Condition)} 
\centerline
{
$F\ =\ \overline{\Int F},$
\qquad\qquad\qquad\qquad\qquad\qquad\qquad\qquad\qquad
}
\medskip
\noindent
and the following  third condition, which is  important for the Dirichlet Problem and regularity,
even when $F$ is linear (i.e., a closed half-space):

\medskip
\noindent
{\bf (Negativity)} 
\centerline
{
$(r,p,A) \ \in\ F \qquad\Rightarrow\qquad (r-s,p, A)\ \in\ F \fa s\geq0.$
\qquad\qquad\qquad
}

\Def{\BB.2}  A closed subset $F\ss \bbj^2 = \bbr\times \rn\times \Symn$ satisfying 
positivity, negativity, and the topological condition will be called a 
{\bf subequation}.  Functions $u\in\Fv(X)$ will be called {\bf $F$-subharmonic
(in the viscosity sense)}.
\medskip

The set $\cp\equiv \{A : A\geq0\}$ is perhaps the most basic example of a subequation.
(It is certainly of  minimal size, up to a translate, given positivity.) 
The smooth $\cp$-subharmonics are defined by $D^2 u\geq0$.  The general $\cp$-subharmonics are 
exactly the classical convex functions (once $u\equiv-\infty$ is excluded).  
Somewhat surprisingly the proof of this is not in the early literature
but is included in [HL$_1$, Prop.2.6] for example. 
Note that this subequation $\cp$ is convex and pure second-order.

\vskip .3in

\centerline{\bf The Distributional  Approach.}
\medskip

The distributions $u\in \cd'(X)$ are continuous linear functionals on the space
$C^{\infty}_{\rm cpt}(X)$ of {\sl distributional test functions}.  For any
linear second-order partial differential operator $L$ with constant coefficients,
$Lu$ is again a distribution.  The notion $u\geq0$ for $u\in \cd'(X)$  is defined by requiring that
$u(\vf)\geq0$ for all $\vf \in C^{\infty}_{\rm cpt}(X)$ with $\vf\geq0$. Thus, the  differential 
inequalities $Lu\geq \l$,  for $\l\in \bbr$, make sense. The pair $L,\l$ defines a half-space
$H(L,\l)$ in $\bbj^2$ by $H(L,\l) \equiv \{ L(r,p,A) \geq\l\}$.

\Def{\BB.3}
Suppose that $F $ is a closed convex subset of $\bbj^2$.
Given $u\in \cd'(X)$, we say that $u\in \Fd(X)$ if 
$$
L u\ \geq\ \l \qquad{\rm whenever\ \ }  H(L,\l) \ {\rm contains\ } F.
\eqno{(\BB.2)}
$$
 
\vskip .3in

\centerline{\bf Comments Concerning the   Approaches.}
\medskip

In two or more variables, the example $F =  \{{\partial^2u\over \partial x_1^2}\geq0 \}$
shows that $\Fv(X)$ and $\Fd(X)$ are in general quite different, and in no sense isomorphic.
For the example of the Laplacian
 $F = \{\D u \equiv \sum_k {\partial^2u\over \partial x_k^2}\geq0\}$ however, 
 $\Fv(X)$ and $\Fd(X)$ are isomorphic. Nevertheless,  as the reader will see in Section 9,
 some of the problems that are
 overcome  in describing an isomorphism between $\Fv(X)$ and $\Fd(X)$   are already
 illustrated by this basic case.
 
 Under suitable hypotheses on $F$,  our main result (Theorem 1.1)  
 obtains an explicit isomorphism between $\Fv(X)$ and $\Fd(X)$. 
 These  hypotheses apply to subequations such as the Laplacian, and   $F \cong \{
{\partial^2u\over \partial x_1^2}\geq0 \} \cup\{ {\partial^2u\over \partial x_2^2}\geq0\}$
on $\bbr^2$ for example.

  \vfill\eject
  

\noindent{\headfont \CC.  The Standard Intersection Theorem for Closed Convex Sets.}
\medskip

Any convex subequation  $F\ss \bbj^2 \equiv \bbr\times \rn\times \Symn$  can be written as an
intersection  of closed affine half-spaces.  The standard way of doing this is as follows.
The dual space
$$
\bbj_2 \ \equiv \ \left (\bbj^2\right)^*
$$
should be viewed as the space of 
 second-order linear partial differential operators $L$ (with constant coefficients)   defined
by
$$
(L\vf)(x)\ =\ \bra a {D^2_x \vf} + \bra b {D_x \vf} + c \vf(x)
\eqno{(\CC.1)}
$$
with $a\in\Symn$, $b\in\rn$, and $c\in \bbr$.  

The {\bf principal symbol} $a=\s(L)$ of a differential operator $L\in \bbj_2$ is defined by
the natural projection $\s :\bbj_2 \to \Symn$ which is dual to the natural inclusion
$\Symn\ss \bbj^2$   (obtained by considering 2-jets of 
functions at a point  which have zero as a critical value at that point).
The constant $c=c(L)$ will be referred to as the {\bf zero$^{\rm th}$-order term}.

Each linear differential operator $L$ determines a closed vector half-space
$H_L$ in ${\bf J^2}$ by requiring $L\vf \geq0$.  
Note that $H_L = \overline{\Int H_L}$ is always true, while:
$$
\eqalign
{
H_L \ \ {\rm satisfies\ (P)} \qquad &\iff\qquad {\rm the\ symbol}\ \ \s(L) =a\geq0   \cr
H_L \ \ {\rm satisfies\ (N)} \qquad &\iff\qquad {\rm the\ zero^{\rm th} order\ term}\ \ c(L) \leq 0   \cr
}
$$
Hence, $H_L$  is a subequation if and only if 
$a=\s(L) \geq0$ and $c(L)\leq0$, in which case it will be referred to as a {\bf linear subequation}.
For each $\l\in\bbr$ the pair $L,\l$ determines the translated half-space subequation
$H(L,\l)$ by requiring
$$
L(r,p,A)\ \equiv \ \bra a {A} + \bra b {p} +cr\ \geq\ \l.
\eqno{(\CC.2)}
$$

Suppose now that $F\ss {\bf J^2}$ is a convex subequation contained in a half-space
$H(L,\l)$ for general $a,b,c$.  
Then it is easy to see that positivity for $F$ implies positivity for $H(L,\l)$, and 
negativity for $F$ implies negativity for $H(L,\l)$.
Thus $H(L,\l)$ is also a subequation.  Said differently, 
$$
F \  \ {\rm is\ a \  subequation\  and \ \ }
F\ss H(L,\l) 
\qquad\Rightarrow\qquad
a = \s(L)\geq0 \quad {\rm and}\quad c(L)\leq0.
\eqno{(\CC.3)}
$$

Although our focus is on convex subequations $F$, which are subsets of the 2-jet space
$\bbj^2 =\bbr\times\rn\times\Symn$ satisfying the additional properties described in Definition  \BB.2,  the elementary constructions
we wish to review hold more generally.  

\vskip .3in
\centerline {\bf The Standard  (Hahn-Banach) Intersection Theorem}
\medskip

We suppose for the moment  that 
 $F$ is any  closed convex subset in a finite dimensional vector space  $V$.
Each closed  affine half-space $H$ can be written parametrically as 
$$
H(w,\l) \ =\ \{v\in V : (w,v)\geq \l\}
\eqno{(\CC.4)}
$$
for some non-zero linear functional $w\in  V^*$ and some $\l\in\bbr$.

\Def{\CC.1}  Given $w\in V^*-\{0\}$ and $\l\in\bbr$, if $F\ss H(w,\l)$,  then
$H(w,\l)$ is called  $F$-{\sl containing} and the linear functional
$w$ is also called $F$-{\sl containing}.
Let $\cc(F)$ denote the set of all  $F$-{\sl containing} linear functionals
(``direction vectors'').
 \medskip

Each proper closed convex set $F$ is the intersection of its $F$-{\sl containing}  half-spaces,
that is
$$
F\ =\ \bigcap   H(w, \l).
\eqno{(\CC.5)}
$$
where $w\in\cc(F)$ and $\l\in\bbr$ is such that $F\ss H(w,\l)$.
This  standard  intersection theorem is a consequence of the geometric form of the Hahn-Banach Theorem (see  Appendix A).
 
It is easy to see that $\cc(F)\cup\{0\}$ is a  convex cone, but it may not be closed.  
Its closure is, of course, a closed convex cone.

\Def{\CC.2. (The Dual Span)}  Let $\cp_+(F)$ denote the closed convex cone $\overline{\cc(F)}$.
The vector space span of $\cc(F)$ (or equivalently $\cp_+(F)$) will be called the 
{\sl dual span of $F$} and denoted by $S_F$.

\medskip

Now back to our case $V=\bbj^2$.  Restating (\CC.5)  for a convex subequation $F$, we have 
$$
F\ =\  \bigcap   H(L, \l).
\eqno{(\CC.5)'}
$$
taken over $F$-containing linear operators $L$ and  $\l\in\bbr$ satisfying $F\ss H(L,\l)$.
By (\CC.3) we also have:
$$
L\in\cc(F) \qquad\Rightarrow\qquad \s(L)\geq 0 \ {\rm and\ } c(L)\leq 0.
\eqno{(\CC.6)}
$$
Unfortunately, unless a convex subequation $F\ss \bbj^2$ is uniformly elliptic, 
 the linear half-space subequations $H \equiv H(L,\l)$ occurring in  (\CC.5)$'$
will typically  be highly degenerate rather than having positive definite symbol.  
The version of (\CC.5)$'$  we need will be proven in the next section (Theorem \DD.2).

If  the principal symbol $\s(L)$ of $L\in \bbj^2$  is positive definite, $H(L,\l)$ is called a 
{\bf  uniformly elliptic linear  subequation}.
One of our main results, Theorem \GGG.2, characterizes those subequations $F$
for which there exists a family $\cf$ of  uniformly elliptic linear  subequations
whose  intersection is $F$.   

The implications of (\CC.5)$'$ can be summarized as follows.

\Remark {\CC.3. (The Two Approaches Revisited)}
Suppose that $F$ is a convex subequation.
\medskip
\noindent
{\bf Viscosity Approach:}   \  Given $u\in \USC(X)$
$$
u\in \Fv(X)
\qquad\iff\qquad
Lu \ {\geq}_{\rm vis} \  \l 
\qquad\forall\, L\in\cc(F) \  {\rm with}\  F\ss H(L,\l).
\eqno{(\CC.7)}
$$
This is a triviality using Definition \BB.1  and  is just as easy to prove in greater
generality.  Given any family of subequations $\{F_\a\}$ (not necessarily convex),
the intersection $F= \bigcap_\a F_\a$ is a subequation, and 
$$
u\in \Fv(X)
\qquad\iff\qquad
u\in \Fv_\a(X) \ \  \forall\  \a
\eqno{(\CC.7)'}
$$

\medskip
\noindent
{\bf  Distributional Approach:}   \  Given $u\in \cd'(X)$, Definition \BB.3 states that
$$
u\in \Fd(X)
\qquad\iff\qquad
Lu \ {\geq}_{\rm dis}\   \l 
\qquad\forall\, L\in\cc(F)\   {\rm with}\  F\ss H(L,\l).
\eqno{(\CC.8)}
$$

\Remark{\CC.4} Even for  the nicest convex  subequations,   the two sets defined above 
can be quite different. For  example, suppose  $F=\cp$
(where $\cc(F)\cup\{0\}=\cp$).  As noted at the end of \S 2, there is no isomorphism between the  set of distrubutions
satisfying $Lv\geq_{\rm dist} 0$ and the set of   u.s.c. functions  satisfying 
$Lu\geq_{\rm visc} 0$  for {\sl all}   $L\in \cc(F)$.
To achieve a bijection we  must restrict  to a subset of $\cc(F)$.

\Lemma{\CC.5}  {\sl
Suppose $\cf \ss\cc(F)$ and the closed convex cone on $\cf$ equals $\cp_+(F) \equiv
\overline{\cc(F)}$.  Then for $u\in\cd'(X)$}
$$
u\in \Fd(X)
\qquad\iff\qquad
Lu \ {\geq}_{\rm dis}\   \l 
\qquad\forall\, L\in\cf \  {\rm with}\  F\ss H(L,\l).
\eqno{(\CC.9)}
$$
The proof is left to the reader.


\vskip .3in

\def\Ir{\Int_{\rm rel}}

\noindent{\headfont \ \DD. Another Intersection Theorem for Closed Convex Sets. }
\medskip

  In this section we suppose that 
 $F$ is an  closed convex proper   subset in a finite dimensional inner product space
$(V, \langle \cdot,\cdot\rangle)$.  
The case where $F$ is unbounded is the case of interest.
Let $\as F$ denote the affine span of $F$. To begin recall that
$$
(1)\ \ \as F \ =\ V\qquad\iff\qquad
(2)\ \ \Int F \ \neq\ \emptyset \qquad\iff\qquad
(3)\ \  F \ =\ \overline{\Int F}.
\eqno{(\DD.1)}
$$
To see that $(1)\Rightarrow (2)$ consider an $n$-simplex obtained from a basis for $V$
contained in $F$ and one other generic point  in $F$.
To see that $(2)\Rightarrow (3)$ note that if $x\in\Int F$ and $y\in \partial F$, then the open segment
joining  $x$ to $y$ must belong to $\Int F$.

 There is no loss in generality in assuming that  the affine span of $F$ is $V$,
or equivalently that $F=\overline{\Int F}$.  For other closed convex subsets $C\ss V$,
let $\Ir C$ denote the interior  of $C$ relative to the affine span of $C$.  Then
we always have that
$$
C\ =\ \overline{\Ir C}.
\eqno{(\DD.2)}
$$
 The set of parameterized   half-spaces $H(w,\l)$ which contain $F$,  
 union with the set $\{0\}\times (-\infty, 0]$,  can be  be written  as   the following subset
 of $V\times\bbr$
$$
C_+(F)\ = \ \{(w,\l) : \bra wv \geq \l \ \forall\, v\in F\}
\eqno{(\DD.3)}
$$
since with $w\neq0$
$$
F\ \ss\ H(w,\l) \qquad\iff\qquad (w,\l)\in C_+(F).
\eqno{(\DD.4)}
$$
This set $C_+(F)$ is obviously a closed convex cone in $V\times \bbr$ with vertex at the origin.

Recall   that by Definition \CC.1   the set of containing direction vectors  $\cc(F)$
for $F$  is simply the projection of $\cc_+(F)$ onto $V$, that is,  
 $$
 \cc(F) \cup \{0\}  \ \equiv\  \pi \{C_+(F)\}
\eqno{(\DD.5)}
$$
  where $\pi : V\times \bbr\to V$. 
 Example A.9 shows that the convex cone  $\pi \{ C_+(F) \}$ may not be closed.

The version of the Hahn-Banach Intersection Theorem needed here involves containing
half-spaces which are in some sense stable.

\Def{\DD.1}  A non-zero vector $w\in V$ is called  {\bf $F$-stable}
if there exists $\l\in\bbr$ such that
$$
(w,\l) \in \Ir C_+(F).
$$
Let $\stab(F)$ denote the set of $F$-stable direction vectors for $F$.
\medskip

\Theorem{\DD.2} 
{\sl Suppose $F$ is a closed convex  proper  subset of $V$. Then $F$ is the 
intersection of the half-spaces    defined by   $F$-stable direction vectors, i.e.,}
$$
F\ =\ \bigcap_{\eqalign
{&w\in \stab(F)  \cr
&F\ss H(w,\l)\cr
}
}
 H(w, \l) 
\eqno{(\DD.6)}
$$
\pf If $v\notin F$, then by (\CC.5) there exists a parameterized containing half-space
$H(w,\l)$ for $F$ which excludes  $v$. That is, 
 $(w,\l)\in C_+(F)$ and  $\bra wv < \l$.  Since this inequality
holds for all $(w',\l')$ in a neighborhood of $(w,\l)$ and $C_+(F) = \overline{\Ir C_+(F)}$,
we may choose $(w',  \l')\in \Ir C_+(F)$ so that  $\bra {w'} v <\l$. 
That is,  $F\ss  H(w', \l')$
and $v\notin H(w', \l')$,  where $w'  \in \stab(F)$ is a relatively stable  direction vector for $F$.
\qed

\Remark{\DD.3} By definition
$$
\stab(F) \ =\ \pi \{\Int_{\rm rel} \cc_+(F)\}.
\eqno{(\DD.7)}
$$
Note  that both $\stab(F)$ and the larger set $\cc(F)$ have the same closure,
denoted $\cp_+(F)$ (see Definition \CC.2).  Hence, they have the same span in $V$,
denoted $S_F$ (the {\sl dual span} of $F$ -- see Definition \CC.2).  Now
$\Int_{\rm rel} \cc_+(F)$ is an open convex cone in $\span \cc_+(F)$.  This is preserved
under the projection $\pi$. That is, $\stab(F)$ is an open convex cone in $S_F$.
Therefore,
$$
\stab(F) \ =\ \Int_{\rm rel} \cp_+(F).
\eqno{(\DD.8)}
$$


\vfill\eject

\noindent{\headfont \EE.  The Edge of  a Convex Set. }
\medskip

Many important examples of convex subequations have a non-trivial 
``edge''  which must be taken into account.

\Def{\EE.1. (The Edge)}
Suppose $F$ is a closed convex subset of $V$.
Fix a point $v_0\in F$.    The  {\bf linearity} or {\bf edge}  
$\Ed(F)$ of $F$ (relative to the point $v_0\in F$)  is
 the set of vectors $v \in V$ such that the line $\ell\equiv \{v_0+tv : t\in\bbr\}$ through
 $v_0$ in the direction $v$ is contained in $F$.  
 \medskip

 Since $F$ is convex, $\Ed(F)$ must be a vector
 subspace of $V$.  The result we need is that:

\Lemma {\EE.2}   {\sl
The vector space  $\Ed(F)$  is independent of the choice of the point $v_0\in F$. 
}

\medskip
This is easy to prove directly by using the convexity of $F$.
This proof is left to the reader.  A second proof follows by 
computing $\Ed(F)^\perp$ and showing this is independent of 
$v_0\in F$.  Recall the dual span $S_F \equiv \span\cc(F)$.

\Prop{\EE.3}  {\sl
$$
\Ed(F)^\perp\ =\  S_F.
$$
This implies that in the decomposition $V = \Ed(F) \oplus S_F$ we have 
$F=\Ed(F)\times F'$ where  $F'$ is  a closed convex subset of 
the dual span $S_F=\span(\stab(F))$ without an edge.
}

\pf
It suffices to show that 
$$
{\rm An \  affine \ line\ } \ell\equiv \{v_0+tv : t\in\bbr\}  \ {\rm   is \ contained \ in \ } F\
\quad\iff\quad 
v\perp \cc(F).
\eqno{(\EE.1)}
$$

To see this note first that the line $\ell$ is contained in a parameterized half-space $H(w,\l)$
if and only if $v\perp w$.
Thus $\ell \ss H(w,\l)$ if and only if $v\perp w$ for all containing direction vectors $w\in \cc(F)$.
Thus, by   (\BB.8)  we have $\ell\ss F$ if and only if $v\perp  \cc(F)$. \qed

\medskip
The following refinement is needed in the proof of a main result, Theorem \GGG.2.

\Lemma{\EE.4} 
{\sl  
Suppose $\bra vw \geq 0,  \ \forall\ w\in\stab(F)$. If $v\notin \Ed(F)$, then $\bra vw>0,  \ \forall\ w\in\stab(F)$. 
}
\pf
Let $\bar v =\pr(v)$ where  $\pr: V\to  \Ed(F)^\perp = \span(\stab(F))$ is orthogonal projection.
Then $\bar v\neq 0$ since $v\notin  \Ed(F)$.  Fix any $w\in \stab(F)$.  Since $\stab(F)$ is open in
$\span(\stab(F))$, there exists $\e>0$ so that $w-\e\bar v \in \stab(F)$. Hence, 
$0\leq \bra v {w-\e\bar v} =\bra vw -\e|\bar v|^2$, which implies that $\bra vw >0$.\qed

\Remark{\EE.5}  
This lemma can be restated in terms of the closed convex cone
 $\cp_+(F)$ and its polar cone, which we denote
 $\cp^+(F) = \cp_+(F)^0$. (See Appendix A for the definition of the polar.)
 $$
 {\rm For}\ \ v\in \cp^+(F), w\in \Ir \cp_+(F) \ \ 
 {\rm one\ has}\ \  \bra vw >0 \ \ {\rm unless}\ \ v\in\Ed(F).
 \eqno{(\EE.2)}
 $$
Note that $\stab(F) = \Ir \cp_+(F)$  (see (\DD.8)).


\vfill\eject

\noindent{\headfont \FF. Subequations Which are Second-Order Complete.}
\medskip

We now  return to the case where  the vector space $V$ is
$$
{\bf J^2}\ =\ \bbr\times\rn\times\Symn,
$$
the space of 2-jets of functions (at the origin) in $\rn$.  Recall that the  dual space
 ${\bf J_2} \equiv( {\bf J^2})^*$, of linear, second-order differential operators at the origin, has a dual splitting
$$
{\bf J_2}\ =\ \bbr\times\rn\times\Symn
$$
and  a natural projection, called the  {\bf principal symbol},
$$
\s : {\bf J_2} \ \arr\  \Symn.
$$

 First, we consider pure
second-order subequations $F\ss\Symn$.  Let $A\bigr|_W$ denote restriction of a 
quadratic form $A\in\Symn$ to a subspace $W\ss\rn$.

\Def{\FF.1}  A subequation  $F\ss \Symn$  {\bf can be defined using fewer of  the  
variables  in} $\rn$  if there exists a  proper subspace $W\ss\rn$ 
and a subset $F'\ss \Sym(W)$ such that 
$$
A\in F\quad \iff\quad A\bigr|_W \in F',
\eqno{(\FF.1)}
$$
Otherwise we say that $F$  {\sl cannot  be defined using fewer of  the  
variables  in} $\rn$.
\medskip

The general case is reduced to the pure second-order case by considering the 
pure second-order subequations associated with $F$, namely the 
fibres of $F$ in $\Symn$:
$$
F_{r,p}\ =\ \{A\in F: (r,p,A)\in F\}.
$$

\Def{\FF.2}  A subequation $F\ss \bbj^2$  is  {\bf second-order complete}
 if $F$ has at least one   fibre $F_{r,p} \ss \Symn$ which cannot  be defined 
 using fewer of  the   variables  in $\rn$.

 \medskip
 
 This definition is clarified by the next result.
 
 \Lemma{\FF.3}  {\sl
 Suppose that $F \ss \bbj^2$  is a convex subequation.  Then $F$ has at least one
 fibre $F_{r,p}$ which cannot  be defined using fewer of  the  
variables  in $\rn$  if and only if 
 every  non-empty fibre $F_{r,p}$ has this same property.
 }
 
 \pf
 Since $\{ A\in\Symn : A\bigr|_W =0\}  = \Sym(W)^\perp$
 (the orthogonal complement of the subspace $\Sym(W)$ in $\Symn$), we see that:
 $$
 \eqalign
 {
 F_{r,p} \ \ & {\rm depends\ only\  on \ the \ variables\  in\ } W \ss\rn   \cr
 & \quad\iff \qquad  F_{r,p} + \Sym(W)^\perp \ \ss\ F_{r,p} \cr
 & \quad\iff \qquad   \Sym(W)^\perp \ \ss\ \Ed(F_{r,p} ).
 }
 \eqno{(\FF.2)}
$$
 If $F_{r,p} \neq 0$, then this condition is equivalent to 
 $$
 \{0\}\times \{0\}\times\Sym(W)^\perp \ \ss\ \Ed(F).
 $$
 
 Now Lemma \EE.2 (with $v_0 = (r,p,A) \in \bbj^2=V$) says that:  
 $ F_{r,p} + \Sym(W)^\perp \ss F_{r,p}$ for one  $ F_{r,p} \neq \emptyset$
 implies the same for all $r,p$. This completes the proof.\qed

 \vfill\eject

\vfill\eject

\noindent{\headfont  \ \GGG.   Characterizing When Stable Means Uniformly Elliptic.}
\medskip

Another equivalent way of saying that $F$ is second-order complete is needed
to finish the proof of our reduction to the linear case (Theorem \AA.2).
This time the proof involves more than the convexity of $F$, it also involves the
positivity condition. Given a unit vector $e\in\rn$, let $P_e:\rn\to\rn$ 
 denote orthogonal projection $P_e(x) = \bra xee$ onto the $e$-line.

\Lemma{\GGG.1}  {\sl
A convex subequation $F$ is second-order complete
\ \ $\iff\ \ {\bf P}_e \equiv (0,0,P_e) \notin \Ed(F)$ for all unit vectors $e\in\rn$.
}
\pf
We prove that: $F$ depends only on the variables in $W\ss \rn$ \ \ $\iff\ \ 
{\bf P}_e\in \Ed(F)$ where $W$ and $e$ are orthogonal.  Since $P_e \in \Sym(W)^\perp$,
 the implication $\Rightarrow$ follows from (\FF.2).
 
 Suppose ${\bf P}_e\in \Ed(F)$.  
 We must show that  $\{0\}\times\{0\}\times \Sym(W)^\perp \ss \Ed(F)$, or that 
 $F_{r,p} + \Sym(W)^\perp \ss   F_{r,p}$ for some non-empty $F_{r,p}$.  Let 
 $M  \equiv \ell +\cp$ where $\ell = \bbr\cdot P_e$ is the line generated by $P_e$.  
 Since ${\bf P}_e \in \Ed(F)$ and $F$ satisifies positivity, 
$F_{r,p} +\overline M \ss F_{r,p} $ where $\overline M$ is the closure of $M$.
The proof is completed by showing that:
$$
\Sym(W)^\perp \ \ss\ \overline M \ \equiv \ \overline{\ell + \cp}.
 \eqno{(\GGG.1)}
$$
Suppose $A \equiv \left(\matrix{s & b \cr b^t & 0}\right) \in \Sym(W)^\perp$ using the blocking
induced by $\rn = \bbr\cdot e \oplus W$.  Let $A_\e \equiv \left(\matrix{s & b \cr b^t & \e}\right)$
for $\e>0$.  Then for $t>>0$ sufficiently large, 
$A_\e +tP_e =  \left(\matrix{s +t & b \cr b^t & \e}\right) \equiv P >0$.  Hence,
$A_\e = -tP_e +P\in \ell + \cp =M$.  Since $A_\e\to A$, this proves that $A\in \overline M$.\qed

\medskip

Note that $\Sym(W)^\perp$ is not contained in $\ell+\cp$ since $A \equiv 
  \left(\matrix{ 0 & b \cr b^t & 0}\right) \in \Sym(W)^\perp$, but 
$  \left(\matrix{ t & b \cr b^t & 0}\right) $  is never $\geq0$.

\Theorem{\GGG.2}
{\sl  Suppose  $F\ss \bbj^2$  is a convex subequation.  Then:
\medskip
\centerline
{
Each $F$-stable linear operator $L$ is uniformly elliptic 
}
\medskip
\noindent
if and only if 
$F$ is second-order complete. }

\pf 
Choose $L\in \stab(F)$.  Assume that $F$ is second-order complete,
or equivalently, by Lemma \GGG.1, that  ${\bf P}_e \notin \Ed(F)$ for each unit vector
$e\in\rn$.  Recall that by (\CC.7) we have 
$$
\bra {{\bf P}_e} L  \ =\ 
\bra {\{0\} \times \{0\}\times P_e}{L}  \ =\ 
\bra {P_e} {\s(L)} \ =\ 
\bra {\s(L) e} e \ \geq\ 0.
$$
Since ${\bf P}_e  \notin \Ed(F)$, we conclude that 
$$
\bra {\s(L)e} e = \bra {{\bf P}_e} L \ >\ 0
$$
by Lemma \EE.4.

Conversely, just assume $L\in\stab(F)$   is uniformly elliptic,  i.e., 
$\bra {{\bf P}_e}L = \bra {\s(L)e} e >0$.
Then ${\bf P}_e$ is not orthogonal to $\stab(F)$ and hence
${\bf P}_e \notin (\span \stab(F))^\perp$ which equals $\Ed(F)$ by 
Proposition \EE.3.\qed
\medskip

This also proves that 
$$
\exists \, L\in \stab(F) \ {\rm with\ } \s(L)>0
\qquad\Rightarrow\qquad
\forall\, L \in\stab(F), \ \s(L)>0.
 \eqno{(\GGG.2)}
$$

\vfill\eject


\noindent{\headfont  \ \HH.   The Main Result for Variable Coefficients. }
\medskip

In this section we derive the main isomorphism  result of the paper. We shall work on an open subset
$X\ss \rn$.  However, the results will be formulated in a way that carries over immediately
to general manifolds.

We denote by $J^2(X)$ the bundle of 2-jets of functions on $X$.  The fibre $J^2_x(X)$ at a point
$x\in X$ is defined to be the germs of smooth functions at $x$ modulo those which vanish
to order three at $x$. Using  the  coordinates on $\rn$ this fibre is naturally
identified with $\bbj^2 = \bbr\times \rn\times\Symn$.  There is a natural short exact
sequence
$$
0\ \ \arr\ \ \Sym(T^*X) \ \ \arr\ \ J^2(X)  \ \harr {\pi}{}\ \  J^1(X)  \ \ \arr\ \ 0
\eqno{(\HH.1)}
$$
where the fibre of $\Sym(T^*X)$ at $x$ consists of the 2-jets of functions with critical
value zero at $x$, and where $J^1(X) = \bbr\oplus T^*X$ is the bundle of 1-jets.

The dual  bundle $J_2(X) = J^2(X)^*$ has a dual short exact sequence 
$$
0\ \ \arr\ \ J_1(X)  \ \ \arr\ \ J_2(X)  \ \harr{\s}{} \  \Sym(TX)  \  \arr \ 0.
\eqno{(\HH.2)}
$$
The  sections of $J_2(X)$ are the variable-coefficient linear second-order
differential operators (possibly degenerate), and  $\s$ is the {\bf principal symbol} map.

There are two basic sub fibre-bundles in $J^2(X)$:
$$
\cp \ \ss\ \Sym(T^*X) \and \bbr_-\ \ss\ J^2(X),
$$
 where the fibre $\cp_x$  is the cone of non-negative quadratic forms in $\Sym(T^*_xX)$
 (germs at x with local minimum value zero at $x$), and the fibre  $(\bbr_-)_x$ is germs at
 $x$ of the non-positive constant functions.  We employ the definition from [HL$_3$].

\Def{\HH.1} By a {\bf subequation} on $X$ we mean a closed subset $F\ss J^2(X)$
with the property that under fibre-wise sum, the Positivity Condition
$$
F+\cp\ \ss\ F
\eqno{(P)}
$$
and the Negativity Condition
$$
F+\bbr_- \ \ss\ F
\eqno{(N)}
$$
and  the Topological Conditions
$$
 (i)\ \ F\ =\ \overline{\Int F}, \qquad (ii)\ \ F_x\ =\ \overline{\Int F_x}, \qquad 
 (iii)\ \ \Int F_x\ =\     (\Int F)_x
 \eqno{(T)} 
 $$
(for each $x\in X$)  hold.

\Def{\HH.2}  A subequation $F\ss J^2(X)$ is {\bf convex} if each fibre $F_x$ is convex.
A subequation  $F\ss J^2(X)$ is {\bf second-order complete}
if each non-empty $\Sym(T^*_xX)$-fibre of $F$ 
cannot be defined using fewer of the variables in $T^*_xX$.
\medskip

\Remark{\HH.3} Recall Lemma \FF.3 which says that for fixed $x\in X$,
 if at least one  $\Sym(T^*_xX)$-fibre of $F$ 
cannot be defined using fewer of the variables in $T^*_xX$, then this holds for all non-empty
$\Sym(T^*_xX)$-fibre of $F$.

\medskip

For any convex subequation, all  of the concepts from the previous
sections, such as $E(F_x)$, $C_+(F_x)$, $\stab(F_x)$, etc. make sense point-wise. 
In particular,  we 
can consider the subset
$$
\stab(F) \ \equiv \ \bigcup_{x\in X} \stab(F_x)\ \ss\ J_2(F)
\eqno{(\HH.3)}
$$

 \Def{\HH.4}  Given a convex subequation $F\ss J^2(X)$, a smooth 
 section $L$ of $\stab(F)$ is said to be an {\bf $F$-stable linear operator}.
\medskip

We first state our main theorem with a mild regularity hypothesis which makes the proof transparent.
Then we show that in most cases of interest  this hypothesis is satisfied.

\Def{\HH.5}  A convex subequation $F$ is called {\bf regular} if for each $x\in X$ and $L_x\in \stab(F_x)$, 
there exists a local $F$-stable linear operator $L$ extending $L_x$.

\medskip

The notion of an  upper semi-continuous $\Fv$-subharmonic function, defined via test functions,
 carries over from the constant
coefficient case to this general one (see [HL$_3$]), as does the notion of a distributional 
$\Fd$-subharmonic function.

\Theorem{\HH.6. (Reduction to the Linear Case)}
{\sl
Suppose that $F\ss J^2(X)$ is a regular convex subequation which is second-order
complete.  Then each $F$-stable linear operator $L$ is uniformly elliptic with $c(L)\leq0$.
Moreover, in both the viscosity and distributional cases
$$
\eqalign
{
&u\  {\sl is\ } F \, {\sl subharmonic \ on\ } X 
\qquad\iff\qquad  \cr
Lu\ \geq\  \l \ \ \ {\sl for\ all\ }  &{ \sl locally\  defined\ } F\, {\sl stable\ linear\ 
 operators\ } L\ {\sl with }\ H(L,\l)\supset F
 }
\eqno{(\HH.4)}
$$
}

\pf  Fix $x\in X$ and for each $L_x \in \stab(F)_x$, let $L$ be the local section of 
$\stab(F)$ extending $L_x$. Define $\cf$ to be the family of all such local sections.
The theorem now follows by  applying  Theorem \GGG.2 and then Theorem \DD.2 combined
with Remark \CC.3.\qed
\medskip

We now take up the question of producing criteria which guarantee that a 
convex subequation is regular.
The following is straightforward, but covers many of the most interesting cases.

\Prop{\HH.7} {\sl
Suppose $F$ is a convex subequation on $X$ with the property that
each point $x\in X$ has a neighborhood $U$ and a local affine trivialization
$$
\Phi: J^2(U) \ \arr\ U\times \bbr^N
$$
of the 2-jet bundle over $U$, such that 
$$
\Phi\left(F\bigr|_U\right)\ =\ U\times \bbf
$$
for a convex subset $\bbf\ss\bbr^N$.  Then $F$ is regular.
}
\pf
In  defining the set $\stab(F)$  in Section \CC, we only used the hypothesis
that $F$ was a convex subset of a vector space $V$.  It is thereby easy to
see that under the local trivialization of the dual bundle $J_2(U)$ induced
by $\Phi$, the set $\stab(F)$ is mapped to a set of the form $U\times {\bf F}$.\qed

\medskip

There is another general criterion for regularity. We say that $F$ {\sl has no edge}
if $\Ed(F_x)=\{0\}$ at each point $x\in X$.

\Lemma {\HH.8}  {\sl 
If $F$ has no edge, then $F$ is regular.}

\pf 
 We will use (\EE.2) to show that for each $x\in X$
$$
\stab(F_x) \ =\ \Int \cp_+(F_x) \ \ss\ \Int \cp_+(F).
\eqno{(\HH.5)}
$$
Fix $x_0 \in X$ and choose a local section $v(x)$ of $\Int F$
on a small closed ball $B$ about $x_0$.  This is always possible 
since,  by (T)(iii), $\Int F_x \ss\Int F$.  
Let $\Sigma \ss J^2(X)$ denote the bundle of unit spheres
centered at  the origin in $J^2(X)$.  For each  $x\in B$ and $r>0$, let  $K_x(r)  \equiv \{v \in \Sigma_x : 
v(x) + r v \in F_x\}$.  
Then the compact sets $K_x(r)$ are fibrewise decreasing as $r\to\infty$,
 so that $K_x\equiv \bigcap_{r>0} K_x(r)$
is a compact subset of $\Sigma_x$.  It follows from the definition in Appendix A that 
the cone on $K_x$ is the asymptotic cone $\overrightarrow{{F_x}}$.  It then follows from Corollary A.6 that
$$
\cp^+(F) \ =\ {\rm Cone}(K) \qquad{\rm  (the\ fibrewise \ cone)}.
\eqno{(\HH.6)}
$$
Therefore by (\EE.2)
$$
w\in \Int\cp_+(F_x) 
\qquad\iff\qquad
\bra wv \ >\ 0 \qquad\forall\, v\in K_x.
\eqno{(\HH.7)}
$$
Since $K$ is compact, if  we are  given $w_0 \in \Int \cp_+(F_{x_0})$,  we
can find a neighborhood $N$ of $w_{x_0}$ covering a (smaller)  ball $B$ about $x_0$ so that 
$$
\inf \left \{ \bra vw : w\in N \ \ {\rm and}\ \ v(x)\in K \ \forall\, x\in B  \right \} \ \geq\ \e
\eqno{(\HH.8)}
$$
for some $\e>0$.  This proves (\HH.5), which easily implies regularity.\qed

\medskip

The following  result generalizes Proposition \HH.7 above.
It is proved by a   straightforward reduction   to the ``no edge'' case in Lemma \HH.8.

\Prop{\HH.9}  {\sl
If $\Ed(F)$ is a smooth sub-bundle of $J^2(X)$, then $F$ is regular.}

\Note{\HH.10}  An important class of subequations which satisfy the hypothesis of 
Proposition \HH.7 (and Prop. \HH.9) are those which are 
locally affinely jet-equivalent to constant coefficient
subequations.  These include many interesting subequations on manifolds.  The reader is
referred to [HL$_3$] for definitions and examples.

\vfill\eject


\noindent{\headfont \II.  The Equivalence of   Various Notions of Subharmonicity for Linear Equations.}
\medskip
  Consider a uniformly elliptic  linear partial differential equation 
 $$
Lu(x) \ =\  {a(x)} \cdot { D^2u(x)} +  {b(x)}  \cdot{ Du(x)} + c(x) u(x) \ =\ \l(x)
$$
where $a, b, c$ and $\l$ are $C^\infty$ on  an open set $X\ss \rn$, and $a>0$ 
is positive definite  and $c\leq 0$ at each point.
The differential inequality $Lu\geq \l$ defines a variable coefficient linear subequation
$H=H(L,\l)$  (the three conditions (P), (N) and (T) are satisfied by $H$), which is of course convex and second-order complete. In this section we outline the proof of Theorem \AA.1 in the uniformly
elliptic linear case where $F$ equals $H\equiv H(L,\l)$.

 In addition to the viscosity and distributional notions of $H$-subharmonicity 
there is a third, more classical notion of $H$-subharmonicity.  For this we must 
define the associated set of $H$-harmonics.
 
We say that $u$ is   {\sl viscosity $H$-harmonic} if   $u$ is $H$-subharmonic and 
$-u$ is $\wt H$-subharmonic, where $\wt H$ is the dual subequation defined by 
$Lu\geq -\l$.
We say that $u$ is   {\sl distributionally $H$- harmonic} if  $Lu=\l$ as a distribution.
In both cases there is a well developed theory of $H$-harmonics.

For example, the $H$-harmonics, both distributional and viscosity, are smooth. 
 This provides the proof  that the two notions of $H$-harmonic are identical. 
 Moreover, there always exists a global $H$-harmonic function which allows
 us to assume $\l\equiv0$ in the proofs of the results we need.
 It  is not as straightforward to make statements relating
the $H$-subharmonics
$\SHV(X)$ and  $\SHW(X)$ since they are composed of different objects.
The bridge is provided by the following  third definition of $H$-subharmonicity.

\Def {\II.1} A function  $u\in \USC(X)$  is {\sl classically $H$-subharmonic} 
if for every compact set $K\ss X$ and every $H$-harmonic function $\vf$ defined 
on a neighborhood of $K$, we have
$$
u\ \leq\ \vf \quad {\rm on} \ \partial K\qquad \Rightarrow
\qquad  u\ \leq\ \vf \quad {\rm on} \ K.
$$
Let $\SHI(X)$ denote the set of these.
\medskip

We always assume that $u$ is not identically $-\infty$ on any connected component of $X$.
We remind the reader that in this section $H$ stands for ``half-space'' and not for ``harmonic''.

In both the viscosity case and the distributional case a great number of results 
have been established.  They essentially include  the following.

\Theorem{\II.2} 
$$
\SHV(X)\ =\ \SHI(X)
$$

\Theorem{\II.3} 
$$
\SHW(X)\ \cong\ \SHI(X)
$$

Note that the first Theorem \II.2 can be stated as an equality since elements of both 
$\SHV(X)$ and $\SHI(X)$ are {\sl a priori} in $\USC(X)$.  By contrast,  the second Theorem \II.3
is not a precise statement until the isomorphism/equivalence is explicitly described.

Theorem \II.3  requires   careful attention in order to apply it to the nonlinear case.
The isomorphism sending $u\in \SHW(X)$ to $U \in \SHI(X)$ is required to  produce the same upper
semi-continuous function $U \in \USC(X)$ independent of the operator $L$.
This is a consequence of the second theorem below.

 We separate out the two directions in Theorem \II.3.  
 Note the parallel with Theorem \AA.1 parts (A) and (B)

\Theorem {\II.3(A)} 
{\sl
If $u\in \SHI(X)$, then $u\in \lloc(X) \ss \cd'(X)$, and as a distribution,  $Lu\geq_{\rm dist}  \l$,   
 that is, $u\in \SHW(X)$.
} 
\medskip

\Theorem {\II.3(B)}
{\sl
If  $u\in \SHW(X)$,  then $u\in \lloc(X)$, and within the $\lloc$-class $u$ of point-wise
defined functions, there exists 
 a unique upper semi-continuous representative    $U \in \SHI(X)$.  It is given by}
$$
U(x) \ \equiv \  \overrightarrow{ {\rm ess} \lim_{y\to x} }  u(y) \ \equiv \ \lim_{r\searrow0} {\rm  ess} \!\! \sup_{B_r(x)} u
\eqno{(\II.1)}
$$
\medskip

The precise statements, Theorem \II.3(A) and Theorem \II.3(B), give meaning to Theorem \II.3.

\Cor{\II.4}  {\sl   If $H\equiv H(L,\l)$ is a uniformly elliptic linear subequation, then both parts 
(A) and (B) of Theorem \AA.1 hold for $F=H$.}

\medskip
\noindent
{\bf  Proof of  Theorem \II.2.}  
We can assume that $\l\equiv 0$. Then the maximum principle applies to $\SHV(X)$.

We first show that $\SHV(X)  \ss \SHI(X)$.
Assume $u\in \SHV(X)$ and $h\in C^\infty$ is $H$-harmonic on a neighborhood of a compact set
$K\ss X$   with  $u\leq h$ on $\partial K$.
Since $\vf$ is a test function for $u$ at a point $x_0$ if an only if $\vf-h$ is  a test function for $u-h$ at $x_0$,
and since $L(\vf-h) = L(\vf) \geq 0$ at $x_0$,
we have $u-h \in \SHV(X)$.  Therefore, the maximum principle applies to $u-h$,
and we have $u\leq h$ on $K$. Hence, $u\in  \SHI(X)$.

Now suppose $u\notin \SHV(X)$.
Then there exists $x_0\in X$ and a test function
$\vf$ for $u$ at $x_0$ with $(L\vf)(x_0) < 0$.  We can assume (cf. [HL$_3$, Prop. A.1]) that $\vf$ is a quadratic  and
$$
\eqalign
{
u-\vf\ &\leq\  - \a |x-x_0|^2 \quad {\rm for} \ |x-x_0|\leq \rho\ \ {\rm and}\cr
&= \ 0\qquad\qquad\qquad {\rm at} \ x_0
}
$$
for  some $\a,\rho>0$.  Set $\psi \equiv -\vf +\e$ where $\e= \a \rho^2$.
  Then $\psi$ is (strictly) $H$-subharmonic on a neighborhood of $x_0$.
Let $h$ denote the solution to the Dirichlet Problem 
for the equation $L(h) = 0$ on  $B\equiv B_{\rho}(x_0)$
with boundary values $\psi$. Since $h$ is the Perron function for $\psi\bigr|_{\partial B}$
and $\psi$ is $L$-subharmonic on $B$, we have $\psi\leq h$ on $\overline B$.
Hence, $-h(x_0) \leq -\psi(x_0)=\vf(x_0) -\e < u(x_0)$.  However,
on $\partial B$ we have $u \leq \vf -\a \rho^2 = -\psi  
=  -h$.  Hence, $u \notin \SHI(X)$.\qed

\medskip

\noindent
{\bf Outline for  Theorem \II.3(A).}  
This theorem is part of classical potential theory, and a proof can be found in 
[HH], which also treats the hypo-elliptic case.
For uniformly elliptic operators $L$ we outline the argument which proves that $u\in\lloc(X)$,
for later use.

Consider  $u\in \SHI(X)$. Fix a ball  $B\ss X$,
and let $P(x,y)$ be the Poisson kernel for the operator $L$ on $B$ (cf. [G]).
Then we claim that for $x\in\Int B$,
$$
u(x) \ \leq\ \int_{\partial B} P(x,y) u(y) d\s(y)
\eqno{(\II.2)}
$$
where $\sigma$ is standard spherical measure.
To see this we first note that for $\vf \in C(\partial B)$, the
 unique solution to the Dirichlet problem for an $L$-harmonic function on $B$ with
boundary values $\vf$ is given by  $h(x) =  \int_{\partial B} P(x,y) \vf(y) d\s(y)$.
Since $u\in \SHI(X)$ we conclude that 
$$
u(x) \ \leq\ \int_{\partial B} P(x,y) \vf(y) d\s(y)
$$
for all $\vf   \in C(\partial B)$  with $u\bigr|_{\partial B}\leq \vf$.
The inequality (\II.2) now follows since $u\bigr|_{\partial B}$ is u.s.c.,  and 
 $u\bigr|_{\partial B} = \inf\{\vf \in C(\partial B): u\leq \vf\}$.

Note that the integral (\II.2)  is well defined (possibly $=-\infty$) since $u$ is bounded above.

Consider a family of concentric balls $B_r(x_0)$ in $X$ for $r_0\leq r\leq r_0+\kappa$
and suppose $x\in B_{r_0}$.  Then for any probability measure $\nu$ on the interval $[r_0, r_0+\kappa]$
we have 
$$
u(x) \ \leq\ \int_{[r_0, r_0+\kappa]} \int_{\partial B_r}  P_{r}(x,y) u(y) d\s(y) \, d\nu(r)
\eqno{(\II.3)}
$$
where $P_r$ denotes the Poisson kernel for the ball $B_r$.
Let $E\ss X$ be the set of points $x$ such that $u$ is $L^1$ in a 
neighborhood of $x$.  Obviously $E$ is open.  Using (\II.3) one  concludes that
if $x\notin E$, then $u\equiv -\infty$ in a neighborhood of $x$ (cf.  [Ho,  Thm. 1.6.9]).  
Hence both $E$ and its complement
are open. Since we assume that $u$ is not $\equiv-\infty$ on any connected component of $X$,
we conclude that $u\in\lloc(X)$.

That  $Lu\geq_{\rm dist} 0$  is exactly Theorem 1 on page 136 of [HH].\qed

\medskip

\noindent
{\bf Proof of  Theorem \II.3(B).}  In a neighborhood of any point $x_0\in X$ the distribution $u \in\SHW(X)$ 
is the sum of an $L$-harmonic function and a Green's potential
$$
v(x) \ =\ \int G(x,y) \mu(y)
\eqno {(\II.4)}
$$
where $\mu\geq 0$ is a non-negative measure with compact support.
Here $G(x,y)$ is the Green's kernel for a ball $B$ about $x_0$.  It suffices to prove
Theorem \II.3(B) for Green's potentials $v$ given by (\II.4).  The fact that $v\in L^1(B)$ 
is a standard consequence of the fact  that $G \in L^1(B \times B)$ 
with singular support  on the diagonal. Since $G(x,y)\leq 0$, (\II.4) defines a point-wise function
$v(x)$ near $x_0$ with values in $[-\infty, 0]$.  By replacing $G(x,y)$ with the continuous kernel
$G_n(x,y)$, defined to be the maximum of $G(x,y)$ and $-n$, the integrals 
$v_n(x) = \int G_n(x,y)\mu(y)$ provide a decreasing sequence of continuous functions converging
to $v$.  Hence, $v$ is upper semi-continuous.  The maximum principle applied to $v-h$ proves
 that $v\in \SHI(X)$.

Finally we prove that if $u\in \lloc(X)$ has a representative $v\in \SHI(X)$, then $v= U$,
the function defined by (\II.1). Since 
$$
\esssup {B_r(x)} u \ =\ \esssup {B_r(x)} v \ \leq \ \sup_{B_r(x)} v,
\eqno{(\II.5)}
$$
and $v$ is upper semi-continuous, it follows that $U(x) \leq v(x)$. 

Note that since $L(-1) = -c\geq0$, the constant function $-1$ is $H$-subharmonic.
Therefore,  $\int P_r(x,y) d\s(y) \leq 1$  for each $r$.  Since $v\in \SHI(X)$, we can  apply (\II.3) to $v$ 
and conclude that 
 $$
 \eqalign
 {
v(x_0) \ &\leq\   {1\over \kappa} \int_{[0, \kappa]} \int_{\partial B_r} P_{r}(x_0,y) v(y) d\sigma(y) \, dr \cr
&\leq \ \left ( \essup{B_\kappa} v\right)   {1\over \kappa}  \int_{[0, \kappa]}  \int_{\partial B_r}  P_{r}(x_0,y)   d\sigma(y) \, dr \cr
& \leq \ \essup{B_\kappa} v \ =\ \essup {B_\kappa} u,  \cr
}
$$
proving that $v(x_0) \leq U(x_0)$.
\qed

\Remark{\II.5}
The construction of $U$ above is quite general and enjoys several nice properties, which we include here.
To any function $u\in \lloc(X)$ we can associate its {\sl essential upper semi-continuous regularization}
$U$ defined by (\II.1).
This regularization $U$ clearly depends only on the $\lloc$-class of $u$.

\Lemma {\II.6} {\sl For any $u\in \lloc(X)$, the function $U$ is upper semi-continuous. 
Furthermore, for any $v\in \USC(X)$ representing the $\lloc$-class  $u$, we have   $U\leq  v$,
and if $x\in X$ is a Lebesgue point for $u$ with value $u(x)$,
then $u(x)\leq U(x)$.}

\pf
To show that $U$ is upper semi-continuous, i.e., $\limsup_{y\to x}  U(y) \leq U(x)$, it 
suffices to show that  
$$
\sup_{B_r(x)} U  \leq \esssup {B_r(x)} u
$$ 
and then let $r\searrow 0$.  
However, if $B_\rho(y) \ss B_r(x)$, then 
$$
U  (y)\ =\ \lim_{\rho\to0} \esssup {B_\rho(y)} u\leq \esssup  {B_r(x)} u.
$$

Letting $r\searrow0$ in (\II.5) proves that $U (x)\leq v(x)$.

For the last assertion  of the lemma suppose that $x$ is a Lebesgue point for $u$ with value $u(x)$, i.e., by definition
$$
\lim_{r\to0} {1\over \left| B_r(x)\right|} \int |u(y)-u(x)| \, dy \ =\ 0,
\qquad{\rm which \ implies} \qquad 
u(x) \ =\ \lim_{r\to0} {1\over \left| B_r(x)\right|} \int u(y) \, dy.
$$
Then $u(x) \leq  \lim_{r\to0} { \rm ess}\sup_{B_r(x)} u \ =\  U(x)$. \qed

\Remark {\II.7} The Definition \II.1 of being ``sub the harmonics'' makes sense
for any subequation $F$,  and the analogue of Theorem \II.2, namely
$$
\Fv(X) \ =\ F^{\rm class}(X),
\eqno{(\II.6)}
$$ 
holds in much greater generality.  In particular, {\sl If  comparison holds
for $F$ and if the Perron functions for the Dirichlet Problem for the dual
subequation $\wt F$ are always $\wt F$-harmonic, then (\II.6) holds.}
Comparison is used to prove 
$\Fv(X) \ss F^{\rm class}(X)$ (the harmonics may not be smooth, so they cannot
necessarily be absorbed into the test function). Our proof that 
$F^{\rm class}(X) \ss \Fv(X)$ for $F=H$ adapts, using that Perron functions 
are $\wt F$-harmonic.

\medskip

 As noted in the introduction,  combining  Theorem \HH.6 with the 
isomorphism Theorems \II.2 and \II.3  proves 
Theorem \AA.1.

\vfill\eject


\centerline{\headfont  \ Appendix A.   Review of the Relevant Convex Geometry. }
\medskip

We have discussed certain aspects of the geometry of unbounded convex sets
which were needed for the proof of the main results, Theorems \DD.2 and \GGG.2.
  In this appendix we give a more complete and
rounded discussion of the relevant convex geometry which should illuminate the previous
abbreviated discussion. This  is   relevant to other aspects
of  the theory of convex subequations. For instance, the notion of boundary convexity needed 
for the Dirichlet problem is defined using the asymptotic interior of $F$ even when
$F$ is not convex.  Our discussion is brief but complete, including proofs.

\vskip .3in

\centerline{\bf The Hahn-Banach Theorem.}
\medskip

The geometrical form is for open sets which we assume to be in a finite dimensional vector space $V$.

\Theorem{A.1}  {\sl
If $X\ss V$ is an open convex set and $z\notin X$, then there exists an affine hyperplane 
$W$ containing $z$ which does not meet $X$.
}
\medskip
For completeness we recall the standard proof.
\pf
Suppose $W$ is an affine plane disjoint from $X$, containing $z$, and of minimal 
codimension $\geq 2$.  Assume $z=0$.  Using the projection map $\pi:V\to V/W$,
we replace $X$ by $\pi(X)$ and $V$ by $V/W$. This reduces to the case where
$W=\{0\}$ and $\dim(V)\geq2$.  
We can assume further that $\dim V=2$ by taking a non-empty slice 
of $X$ with a 2-dimensional subspace.
Now let $\wt X$ denote the radial projection of
$X$ onto  onto the unit circle in $V$. Then $\wt X$ is open, connected and 
contains no antipodal pairs.  Hence, there must be a line through the origin disjoint
from $X$, which provides a contradiction to codim$(W)\geq2$.
\qed
\medskip

Note that the openness of  linear and radial projections is used in the proof.

There is more than one standard result of this type when the convex set is
closed instead of open, but they are all easy consequences of the open case.
We shall state them using the set $\ch_F$ of (closed)  containing half-spaces  for $F$ (see \S \CC).

 \Theorem {A.2} {\sl
 Suppose that $F$ is a closed convex subset of $V$.
 \medskip
 
 (1) If $z\in\partial F$, then $\exists \, H\in \ch_F$ with $z\in \partial H$.\medskip
 
 (2) If $z\notin F$, then
 \smallskip
 
 \qquad (a) $\exists\, H\in\ch_F$ with $z\notin H$ and $F\cap\partial H \neq \emptyset$.

 \smallskip
 
 \qquad (b) $\exists\, H\in\ch_F$ with $z\notin H$ satisfying $F\ss \Int H$.
}

\medskip

Given a containing half-space $H$ for $F$ ($H\in\ch_F$), if $\partial H \cap F \neq \emptyset$,
then $H$ is {\bf supporting} for $F$, while if $F\ss \Int H$, then $H$ is {\bf strictly containing} for 
$F$.  Hence, the $H$'s in (1) and (2a) are supporting, while the $H$ in (2b) is strictly containing.

\pf First assume that $F=\overline{\Int F}$. Then (1) is immediate from Theorem A.1 
since $z\notin X\equiv  \Int F$. For (2) choose $x\in\Int F$ 
and consider the line segment joining $x$ to $z$.
It has a unique interior point $y$ on $\partial F$.  Now using Theorem A.1 choose a hyperplane $W$
disjoint from $\Int F$ and containing $y$.  Taking $H$ to be the side of $W$ containing $F$
proves 2(a), because this $H$ is supporting (i.e., $y\in F\cap \partial H$) and $x\in \Int H, y\in \partial H \Rightarrow z\notin H$.
Translating $H$ in the direction $z-y$ a small amount produces a different $H$ satisfying (b).

 More generally, if $z\in {\rm Aspan}(F)$, the affine span of $F$, 
 then by  (\DD.1) we can assume that $\span (F) =V$.
Recall that (\DD.1) says that 
$$
{\rm Aspan}(F) \ =\ V 
\qquad\iff\qquad
\Int F \ \neq\ \emptyset 
\qquad\iff\qquad
F\ =\ \overline{\Int F}.
$$

This proves part (1) in general, and part (2) if $z\in {\rm Aspan}(F)$.

To prove part (2) for $z\notin {\rm Aspan}(F)$, only consider half-spaces $H$ with boundary
$W$ parallel to $ {\rm Aspan}(F)$.  The proof of (a) and (b) is straighforward using these 
parallel half-spaces.
\qed

\vskip .3in

\centerline{\bf The Bipolar Theorem.}
\medskip

We define the {\bf polar} of an arbitrary set $X\ss V$ to be 
$$
X^0\ \equiv \ \{w\in V^* :  (w,x)\geq -1\ \forall\,x\in X\}.
\eqno{(A.1)}
$$
Note that $X^0$ is a closed convex set containing 0. Moreover, if
$F\equiv$ the closed convex hull of $F\cup\{0\}$, then
$$X^0=F^0.$$
In the definition (A.1) of the polar, we chose to use $(w,x)\geq -1$ rather than 
$(w,x)\leq 1$.  This ensures that 
If $X$ is a cone with vertex at the origin, then
$$
X^0\ \equiv \ \{w\in V^* :  (w,x)\geq0\ \forall\,x\in X\}.
\eqno{(A.2)}
$$

\Theorem{A.3}  {\sl
If $F$ is a closed convex set containing the origin in $V$, then}
$$
(F^0)^0\ =\ F
$$
\pf
That $F\ss (F^0)^0$ is obvious.  If $z\notin F$, then by Theorem A.2, part 2(b), there exist
$w\in V^*, w\neq0$ and $\l\in\bbr$ such that $F\ss \{x:(w,x)>\l\}$ but $(w,z)<\l$.
Now $0\in F\Rightarrow \l<0$. 
Replacing $w$ by $\overline w = w/|\l|$, we have $\overline w \in F^0$.
Since $(\overline w,z)<-1$ we have $z\notin (F^0)^0$.\qed

\vskip .3in

\centerline{\bf Intersection Theorems.}
\medskip

We recall the notation of  Section \CC.  (See the discussion preceding (\CC.6).)
Let $\cc(F)$ denote the set of containing linear functionals (direction vectors)
for $F$.  That is, $w\in\cc(F) \ss V^*-\{0\}$ if 
$$
F\ \ss\ H(w,\l) \ \equiv\ \{v\in V : (w,v)\geq\l\} \qquad{\rm for\ some\ \ } \l\in\bbr.
$$
For each $w\in\cc(F)$, set $\l_w = \sup\{\overline \l : H(w,\overline \l) \supset F\}$
so that $H(w,\l_w)$ is the smallest closed half-space containing $F$.  By Theorem A.2
$$
F\ =
 \bigcap_{w\in\cc(F)}  H(w,\l_w).
\eqno{(A.3)}
$$

In fact, this remains true for the  smaller set of supporting half-spaces.
We say that a direction vector $w$ is  {\bf supporting} if $F\cap \partial H(w,\l_w)\neq \emptyset$. 
Let $\spt(F)$ denote the subset of $\cc(F)$ consisting of  these  supporting direction vectors $w$ for 
$F$. Then also by Theorem A.2
$$
F\ =
 \bigcap_{w\in\spt(F)}  H(w,\l_w).
\eqno{(A.4)}
$$
This is the standard intersection result (\CC.6).

Note that if $F$ is a bounded set, then all non-zero $w\in V^*$ are supporting for
$F$, and $\spt(F) = \cc(F) = V^*-\{0\}$, so there is no difference between (A.3) and (A.4)
in this case.

\vskip .3in

\centerline{\bf A More Complete Picture.}
\medskip

We fix an (unbounded) closed convex  proper subset $F$ 
 in a finite dimensional inner product space $V$. We assume $\Int F\neq \emptyset$.
We shall associate to $F$ several  closed convex cones.  They come in polar pairs:
first $C^+(F)$ and $C_+(F)$,  then $\cp^+(F)$ and $\cp_+(F)$, and finally the polar subspaces
$ \Ed(F)$ and $S_F$ (the dual span of $F$).

We start with the cone  in $V\times \bbr$ which describes 
the parameterized containing half-spaces for $F$:
$$
C_+(F)\   =  \  \{(w,\l) : F\ss H(w,\l)\}\ \ss\ V\times \bbr,
\eqno{(A.5)}
$$
but also includes $\{0\}\times [-\infty,0]$
Equivalently, this can be written as
$$
C_+(F)\   =  \  \{(w,\l) :  \bra wv -\l\geq0 \ \forall\, v\in F\}.
\eqno{(A.5)'}
$$
By definition, this is the polar of the set $F\times \{-1\}$ in the vector space $V\times \bbr$.
It is convenient at this point to define $C^+(F)$ to be the polar of $C_+(F)$.
The Biploar Theorem together with (A.5)$'$ states that
$$
C^+(F)\   =  \  \overline{\Cn(F\times \{-1\})}.
\eqno{(A.6)}
$$

Our first polar pair can be further  analyzed as follows.  Note that the convex cone
$\Cn((F\times \{-1\})$ is closed when considered as a subset of the open half-space
in $V\times \bbr$ defined by $\l<0$.  We define $\cp^+(F)$ to be the remaining part
of $\overline{\Cn(F\times \{-1\})}$, namely, 
its intersection with the boundary hyperplane $V\times \{0\}$. In other words,
$$
\cp^+(F)\   \equiv  \  \{v\in V : (v,0) \in C^+(F)\}.
\eqno{(A.7)}
$$
Then
$$
C^+(F)\   =  \   \Cn(F\times \{-1\}) \cup (\cp^+(F) \times \{0\})
\eqno{(A.8)}
$$
is a disjoint union except for the origin.

Of course, the intersection of $\cc^+(F)$ with the hyperplane $\{\l=-1\}$ is $F\times \{-1\}$.
However, note that the intersection of $\cc_+(F)$ with   $\{\l=-1\}$ is $F^0\times \{-1\}$.

As with $C^+(F)$ and $C_+(F)$, it is convenient to define $\cp_+(F)$
to be the polar of $\cp^+(F)$.  It is trivial to check that
$$
A^0 \cap B^0\ =\ (A+B)^0
\eqno{(A.9)}
$$
for two closed convex cones $A$ and $B$.  Hence,
$$
 (A\cap B)^0   \ =\ \overline{ (A^0  + B^0)}
\eqno{(A.9)'}
$$
follows from the Bipolar Theorem.
Since $\cp^+(F) \times \{0\}$ is the intersection of $C^+(F)$ with $V\times \{0\}$, its polar $\cp_+(F) \times \bbr$ 
satisfies
 $$
\cp_+(F) \times \bbr  \ =\ \overline{C_+(F) + (\{0\}\times \bbr)}.
\eqno{(A.10)}
$$
Let $\pi : V\times \bbr \to V$ denote projection. Then (A.10) can be rewritten as
 $$
\cp_+(F) \ =\ \overline{\pi\{C_+(F)\}},
\eqno{(A.10)'}
$$
the closure of $\cc(F) \equiv \pi\{C_+(F)\}$, the set of {\sl containing direction vectors for $F$}.
This was essentially taken  as the definition of $\cp_+(F)$ in Section \EE.

\Remark{A.4}   In general, for a closed convex cone $C$ the projection $\pi (C)$
may not be closed (see Example A.9 below), but $\pi ( \Int C)$ is always open and has the same closure as  $\pi (C)$.

\medskip

The {\sl stable direction vectors for $F$}
 $$
\stab(F) \ \equiv\ \pi\{\Ir C_+(F)\},
\eqno{(A.11)}
$$
form an open convex cone with vertex (missing) at the origin in the dual span
 $S_F\equiv \span \stab(F)$.  Also note that
 $$
\stab(F) \ \ss \ \cc(F)\ \ss\ \cp_+(F)
\and
\stab(F) \ = \ \Ir \cp_+(F).
\eqno{(A.12)}
$$

Since $\stab(F)$ is the relative interior of $\cp_+(F)$, all three sets
$\stab(F) \ss \cc(F)\ss \cp_+(F)$ have the same vector space span $S_F$.
Our final polar pair is  $\Ed(F)$ and $S_F$.  See Section \EE \ (Prop. \EE.3) for one proof that 
they form a polar pair. Alternatively, note that 
$\Ed(F)^\perp=\span \cp^+(F)$ and $S_F = \span\cp_+(F)$.

\vskip .3in
\centerline{\bf The Asymptotic Cone at Infinity  and the Monotonicity Set for $F$.}\medskip

Our discussion of the asymptotic cone of $F$ is parallel to the discussion of the edge
in Section \EE. To begin we fix a point $v_0\in F$ and define the {\bf asymptotic cone } $\Fa$
for $F$ (at $v_0$) to be the set of vectors $v\in V$ such that the ray
$R=\{v_0+tv : t\geq0\}$ is contained in $F$. Since $F$ is a closed convex set,
it is easy to see that $\Fa$ is a closed convex cone with the vertex at 
the origin.

\Lemma{A.5}
{\sl
Given $v_0\in F$, a ray $R=\{v_0+tv : t\geq0\}$ is contained in 
$F$ if and only if $\bra vw\geq0 \ \forall\, w\in \cc(F)$.
}
\pf By the (Hahn-Banach) Intersection Theorem (A.4)
 it suffices to show that for each parameterized half-space
$H(w,\l)$ containing $F$ we have
$$
R\ \ss\ H(w,\l) \qquad\iff\qquad  \bra vw\ \geq\ 0.
\eqno{(A.13)}
$$
However, since $a\in F$, we have  that:\smallskip

\centerline{
$R\ss H(w,\l) \iff 0\leq \bra  {a+tv}w = \bra aw +t\bra vw \geq \l \iff \bra vw\geq0$.\quad\mathqed}

\Cor{A.6} {\sl Lemma A.5 can be restated as}
$$
\Fa \ =\ \cp^+(F).
$$
{\sl In particular, $\Fa$ is independent of the choice of $v_0\in F$.
}
\medskip

A vector $v\in V$ is called a {\bf monotonicity vector} for $F$ if
$$
F+v\ \ss\ F.
$$
The set $M(F)$ of such vectors is called the {\bf monotonicity set}  for $F$.
It is easy to show that $M(F)$ is a closed convex cone.

 \Lemma{A.7}
 $$
M(F) \ =\  \cp^+(F).
$$
 \pf
 Choose $v\in \cp^+(F)$. Then by Corollary A.6 we have $v_0+v \in F$. Hence, $\cp^+(F)\ss M(F)$.
 
 Given  $v\in M(F)$ and $v_0\in F$, it follows by induction that $v_0+kv\in F$ for each $k\in \bbz^+$.
 By the convexity of $F$ this implies that the full ray $\{v_0+tv : t\geq0\}\ss F$.\qed

 \medskip
 
Finally we prove that $\stab(F) \ss\spt(F)$ and add an example where 
$\cc(F) \equiv\pi C_+(F)$ is not closed. (See the discussion prior to (A.4) for the
definition of the set $\spt(F)$ of supporting direction vectors.)
 
 \Lemma{A.8}
 {\sl
 Each stable direction vector for $F$ is a supporting direction vector for $F$.
 }
 
 \pf Suppose $w\in\cc(F)$ but $w\notin \spt(F)$, i.e., $H\equiv H(w, \l_w)$ 
 is not a supporting half-space for $F$.
 We will show that $w\notin \stab(F)$.
 Note that by definition $\partial H\cap F=\emptyset$ but $\dist(\partial H,F)=0$.
Hence we can choose a sequence of points $\{z_k\}$ in $F$, going to $\infty$, with
$\dist(z_k, \partial H)\to 0$.  We can assume  that $z_k/|z_k| \to v$.  Then for any point $a\in V$ the line 
segments $[a,z_k]$ converge to the ray $\{a+tv : t\geq0\}$.  Taking $a\in \partial H$ and then $a\in F$
we obtain parallel rays in the direction $v$, one contained in $\partial H$ and the other contained
in $F$.  Hence $v\perp w$ and $v\in \Fa$ proving that $w\in\partial \cp_+(F)$, i.e.,  $w\notin \stab(F)
= \Ir \cp_+(F)$.\qed

\Ex{A.9} This  elementary example  is easy to visualize.
Let $F \equiv \{y-{x^2\over 2} \geq0 \}\ss\bbr^2$.  Then $C^+(F)\ss\bbr^3$ is defined by
$-2\l y - x^2 \geq0$ and $\l\leq0$.  This circular convex cone is uniquely determined 
by the facts that $e\equiv (0,1,-1)$ generates the center ray while $e_2 \equiv (0,1,0)$ and 
$-e_3\equiv (0,0,-1)$ generate boundary rays. This is self-polar, that is, $C_+(F) =C^+(F)$.
One can compute that $\cp^+(F)$ is the ray through $e_1 \equiv (0,1) \in\bbr^2$
by intersecting $C^+(F)$ with $\bbr^2\times \{0\}$, or by noting that this is the asymptotic cone at infinity
(or monotonicity set) for $F$. Its polar $\cp_+(F)$ is the closed upper half-plane
$\{y\geq0\}$.  The set $\stab(F) \equiv \pi(\Int C_+(F))$ of stable directions vectors for $F$ is easily seen
to be the open upper half-plane $\{y>0\}$.
Also the set of containing direction vectors $\cc(F)$ equals the set
$\stab(F)  = \{y>0\}$.  Finally note that $\pi\{C_+(F)\} = \cc(F) \cup \{0\}$, providing an example
where  $\pi\{C_+(F)\}$ is not closed.

\vskip .3in



\centerline{\bf References}

\vskip .2in



\noindent
\item{[C]}   M. G. Crandall,  {\sl  Viscosity solutions: a primer},  
pp. 1-43 in ``Viscosity Solutions and Applications''  Ed.'s Dolcetta and Lions, 
SLNM {\bf 1660}, Springer Press, New York, 1997.

 \smallskip

\noindent
\item{[CIL]}   M. G. Crandall, H. Ishii and P. L. Lions {\sl
User's guide to viscosity solutions of second order partial differential equations},  
Bull. Amer. Math. Soc. (N. S.) {\bf 27} (1992), 1-67.

 \smallskip

 \noindent
\item {[G]} P. Garabedian, { Partial Differential Equations},    J. Wiley and Sons,  New York, 1964.

\smallskip

\item {[HL$_{1}$]}   F. R. Harvey and H. B. Lawson, Jr.,  {\sl  Dirichlet duality and the non-linear Dirichlet problem},    Comm. on Pure and Applied Math. {\bf 62} (2009), 396-443. ArXiv:math.0710.3991

\smallskip

\item {[HL$_{2}$]}  \ \----------,    {\sl  Plurisubharmonicity in a general geometric context},  Geometry and Analysis {\bf 1} (2010), 363-401. ArXiv:0804.1316.

\smallskip

\item {[HL$_{3}$]}  \ \----------,   {\sl  Dirichlet duality and the nonlinear Dirichlet problem
on Riemannian manifolds},  J. Diff. Geom. {\bf 88} (2011), 395-482.   ArXiv:0912.5220.
\smallskip

\item {[[HL$_{4}$]}  \ \----------, {\sl  The restriction theorem for fully nonlinear subequations}, 
ArXiv:1101.4850.

\smallskip

\item {[HL$_{5}$]}  \ \----------,   {\sl  Potential Theory on almost complex manifolds}, 
Ann. Inst. Fourier (to appear).  ArXiv:1107.2584.
\smallskip

\item {[HL$_{6}$]}  \ \----------,   {\sl  Existence, uniqueness and removable singularities for nonlinear
partial differential equations in geometry}, 
Surveys in Geometry(to appear).  
\smallskip

\item {[HH]}  M. Herv\'e and R.M. Herv\'e., {\sl  Les fonctions surharmoniques  dans l'axiomatique de
M. Brelot associ\'ees \`a un op\'erateur elliptique d\'eg\'en\'er\'e},    Annals de l'institut Fourier,  {\bf 22}, no. 2 (1972), 131-145.

\smallskip

 \noindent
\item{[Ho]}
L. H\"ormander,  An introduction to complex analysis in several variables,  Third edition. North-Holland Mathematical Library, 7. North-Holland Publishing Co., Amsterdam, 1990. 
 


   \noindent
\item{[I]}    H. Ishii,  {\sl On the equivalence of 
two notions of weak solutions, 
viscosity solutions and distribution solutions},
Funkcial. Ekvac. 38 (1995), no. 1, 101Ð120.

   \smallskip

   \noindent
\item{[K]}    N. V. Krylov,    {\sl  On the general notion of fully nonlinear second-order elliptic equations},    Trans. Amer. Math. Soc. (3)
 {\bf  347}  (1979), 30-34.

\smallskip

\item {[P]}  N. Pali, {\sl Fonctions plurisousharmoniques et courants positifs de type (1,1)
sur une vari\'et\'e presque complexe}, 
Manuscripta Math.  {\bf 118} (2005), no. 3, 311-337.

\end